\documentclass[11pt]{article}
\usepackage{mathrsfs}
\usepackage{amssymb}
\usepackage{amsmath}
\usepackage{amsthm}
\usepackage{amsfonts}
\usepackage{color,xcolor}
\usepackage{graphicx}
\usepackage{epstopdf}
\usepackage{geometry}
\usepackage{manfnt}
\usepackage[pagewise]{lineno}
\usepackage[colorlinks,citecolor=blue,urlcolor=blue]{hyperref}
\newtheorem{theo}{\textbf{\ \ \quad Theorem}}[section]

\newtheorem{remark}{\textbf{\ \ \quad Remark}}[section]

\newtheorem{defi}{\textbf{\ \ \quad Definition}}[section]

\newcommand{\be}{\begin{equation}}
\newcommand{\ee}{\end{equation}}
\newcommand\bes{\begin{eqnarray}}
\newcommand\ees{\end{eqnarray}}
\newcommand{\bess}{\begin{eqnarray*}}
\newcommand{\eess}{\end{eqnarray*}}

{\setlength\arraycolsep{2pt}}

\title{Geometry property of stochastic parabolic equations and Pr\'{e}kopa-Leindler inequality on bounded domains}

\author{Guangying Lv \footnote{Corresponding author. gylvmaths@nuist.edu.cn}\\
\\    {\small \it  School of Mathematics and Statistics,
Nanjing University of Information}\\
{\small \it  Science and Technology, Nanjing 210044, China}\\
}
\date{}

\begin{document}
\maketitle
\begin{abstract}
We study the preservation of spatial concavity and log-concavity for stochastic parabolic equations posed on bounded convex domains. For the linear noise, by using Minkowski convolutions, we establish the convexity preservation of solution almost surely. For the nonlinear noise, by using It\^{o} formula, the impact of noise is given. Moreover, we establish the relationship between the weighted (stochastic) Pr\'{e}kopa-Leindler inequality and properties of (stochastic) parabolic equation on bounded domain, which is new even for deterministic case.

\medskip
{\bf Keywords}: Concavity; Log-concavity; Stochastic parabolic equations; Pr\'{e}kopa-Leindler inequality; Minkowski convolution

\medskip
AMS subject classifications (2020): 35B30; 35K10; 60H15

\end{abstract}

\baselineskip=15pt

\section{Introduction}\label{sec1}\setcounter{equation}{0}

Geometric properties of functions and solutions to partial differential equations
are often intimately linked to analytic regularity.
For instance, Taylor's formula implies that
\[
f\in C^2,\quad f''(x)\ge 0 \ \ \forall x\in\mathbb{R}
\Longrightarrow
f(\lambda x+(1-\lambda)y)\le \lambda f(x)+(1-\lambda)f(y),
\]
for all $x,y\in\mathbb{R}$ and $\lambda\in(0,1)$.
Here the second derivative is a local analytic quantity, whereas convexity is a global
geometric property.
The converse implication, however, generally fails: a convex function need not be twice
differentiable.
This basic observation motivates the study of geometric properties of solutions to
evolution equations, where regularity may be limited or degenerate.
In this paper, we investigate the preservation of concavity and log-concavity for
solutions to stochastic parabolic equations posed on bounded domains.

Concavity and log-concavity play a central role in the qualitative analysis of partial
differential equations.
They are closely related to stability, uniqueness and monotonicity phenomena, and arise
naturally in problems from potential theory, heat propagation and nonlinear diffusion.
A fundamental question in this direction is whether geometric properties of the initial
datum are preserved along the evolution of the equation.

In the deterministic setting, concavity properties of elliptic and parabolic equations
have been extensively studied.
The quasi-concavity of equilibrium potentials for convex bodies was established in
\cite{G1955,G1957,L1977}, while probabilistic methods were employed by Borell
\cite{B1982} to analyze concavity for heat conduction problems.
The systematic theory of log-concave functions was developed by Brascamp and Lieb
\cite{BL1976}.
A general framework for concavity-preserving properties of partial differential
equations was initiated by Korevaar in \cite{K198301,K198302}; see also
\cite{GK1999,K2001,K2005} and the references therein.
Further related results for parabolic and degenerate diffusion equations can be found in
\cite{A1985,CW2021,CW202101,CH2013,CMS2019,IP2008,IP20102,IST2022}, as well as in the
context of porous medium equations
\cite{B1987,CW2024,DHL2001,IP20101,K2003}.

For parabolic equations on bounded domains, Ishige and Salani \cite{IS2014} studied
space-time concavity for Dirichlet boundary value problems under suitable structural
assumptions on the nonlinear terms.
Their approach relies on viscosity solutions combined with a weak comparison principle.
Subsequent extensions were obtained in \cite{INS2016}, but the analysis remains restricted
to power-type concavity indices $p\in(0,1)$.
A different approach was proposed in \cite{INS2020}, where spatial concavity for
parabolic systems was established via an auxiliary parabolic system and a refined
comparison principle.
While effective for spatial concavity, this method does not cover the borderline
log-concave case $p=0$.

The presence of stochastic perturbations introduces substantial additional difficulties.
Random forcing typically destroys pointwise order structures and prevents the direct
application of classical maximum principles.
In particular, for stochastic parabolic equations driven by Brownian noise, one cannot
expect simultaneous preservation of concavity in both space and time.
This raises the natural question of whether spatial concavity or log-concavity can still
be preserved for stochastic parabolic equations on bounded domains, and under what
structural conditions such properties remain valid.

The first objective of the present paper is to address this question.
We consider stochastic parabolic equations of the form
\begin{equation}\label{1.4}
\left\{
\begin{array}{llll}
du(t,x)
=
\bigl[\Delta u(t,x)+f(t,x,u(t,x))\bigr]dt
+\sigma(t,x,u(t,x))\,dB(t),
& x\in D,\ t>0,\\
u(t,x)=0,
& x\in\partial D,\ t>0,\\
u(0,x)=u_0(x),
& x\in D,
\end{array}
\right.
\end{equation}
where $D\subset\mathbb{R}^n$ is a bounded convex domain and $B(t)$ is a one-dimensional
standard Brownian motion.
Our aim is to identify sufficient conditions under which spatial concavity or
log-concavity of the initial datum is preserved along the stochastic evolution, and to
clarify the role played by the noise structure and boundary conditions.

Our analysis distinguishes between linear and nonlinear noise.
In the linear case,  we develop a viscosity-solution framework combined with a
comparison principle to establish almost sure preservation of spatial concavity under
appropriate structural assumptions on the reaction term.
This approach allows us to treat the borderline case $p=0$, which appears to be new even
for deterministic parabolic equations on bounded domains. For nonlinear noise, the F-concavity is studied. By using It\^{o} formula, we obtain the impact of noise. More precisely, the relationship between the reaction term and diffusion is considered to make the solution have the same property as the initial data.
An example is given to illustrate our main result.

The second objective of the present paper is to exploit the connection between the generalized Pr\'{e}kopa-Leindler inequality and properties of parabolic equation on bounded domain. Notice that Pr\'{e}kopa-Leindler inequality is a special case of Borell-Brascamp-Lieb inequality, and we will study it in Section 4. Recently, Ishige et al. \cite{ILS2025} studied the relationship between Borell-Brascamp-Lieb inequality and properties of diffusion equations of porous medium type. \cite{KM2017} considered the Brascamp-Lieb-type inequalities on weighted Riemannian manifolds with boundary. There is only a small amount of work on bounded domain research about the Borell-Brascamp-Lieb inequality. In this paper, inspired by \cite{ILS2025}, the relationship between a weighted Pr\'{e}kopa-Leindler inequality and properties of parabolic equation on bounded domain will be given, which is new even for deterministic case.

\smallskip
Our contributions can be summarized as follows.
\begin{itemize}
\item[(i)] We establish sufficient conditions for the preservation of concavity for the
Dirichlet problem associated with \eqref{1.4}.

\item[(ii)] The impact of noise on concavity is given.

\item[(iii)] We establish the relationship between a weighted (stochastic) Pr\'{e}kopa-Leindler inequality and properties of (stochastic) parabolic equation on bounded domain.
\end{itemize}

The rest of the paper is organized as follows.
Section~2 is devoted to the linear noise.
In Section~3, we investigate the impact of noise on concavity of solutions. In Section~4, we establish the relationship between the generalized Pr\'{e}kopa-Leindler inequality and properties of parabolic equation.

\section{Linear Multiplicative noise}\label{Sec2}
\setcounter{equation}{0}
Inspired by \cite{IS2014,INS2020}, in this section we investigate the Dirichlet boundary
value problem for stochastic parabolic equations of the form
\bes
\left\{
\begin{array}{llll}
d u(t,x)
=
\bigl[\Delta u(t,x)+f(t,x,u(t,x))\bigr]dt
+\sigma(t,x,u(t,x))dB(t),
& x\in D,\ t>0,\\
u(t,x)=0,
& x\in\partial D,\ t>0,\\
u(0,x)=u_0(x),
& x\in D,
\end{array}
\right.
\label{2.1}
\ees
where $D\subset\mathbb{R}^n$ is a bounded convex domain and $B(t)$ is a
one-dimensional standard Brownian motion. If the noise depends on spatial vector, then it is impossible to consider the concavity of spatial vector in the sense of almost surely.

In particular, we focus on noise of the form $\sigma(t,x,u)=\sigma(u)$.
It was shown in \cite{KY2020} that, if $f\in L^r([0,T]\times D)$ and $\sigma$ satisfies a
Lipschitz condition, then the solution of \eqref{2.1} remains positive almost surely; see
\cite[Theorem~4.4]{KY2020}. Let $(\Omega,\mathcal{F},\mathbb{P})$ be a complete probability space with filtration $(\mathcal{F}_t)_{t\in[0,T]}$ generated by $B(t)$. The initial value $u_0$ is a $\mathcal{F}_0$-measurable variable in $L^2(D)$. A predictable $L^2(D)$-valued stochastic process $\{u(t) : t\in [0, T ]\}$ is called a weak solution of \eqref{2.1} if for any $v\in H^1(D)$ and almost every $t\in[0,T]$, the following equality holds
  \bess
(u(t),v)_{L^2(D)}=(u_0,v)+\int_0^t\left[(-u(s),\Delta v)+(f(s,\cdot,u(s)),v)\right]ds+\int_0^t(\sigma(s,\cdot,u(s)),v)dB(s).
  \eess
We say $u$ is a strong solution of \eqref{2.1} if $u\in L^\infty(0,T;L^2(D))\cap L^2(0,T;H^1(D))$ with $H^1(D):=\{u\in L^2(D),\nabla u\in L^2(D)\}$ and $\|u\|^2_{H^1(D)}:=\|u\|^2_{L^2(D)}+\|\nabla u\|^2_{L^2(D)}$.
In this subsection, we investigate the preservation of spatial log-concavity for
solutions to \eqref{2.1}. Due to the technical reason, we only consider the case $\sigma(t,x,u)=\sigma u$, where $\sigma$ is a constant.
We first introduce the following definitions.

\begin{defi}\label{d2.1}
Let $D$ be a convex domain in $\mathbb{R}^n$. A nonnegative function $f$ on $D$ is said to
be $p$-concave in $D$ if
\[
f\bigl((1-\lambda)x+\lambda y\bigr)\geq M_p\bigl(f(x),f(y);\lambda\bigr)
\]
for all $x,y\in D$ and $\lambda\in(0,1)$, where $M_p(a,b;\lambda)$ denotes the $p$-mean of
two nonnegative numbers $a$ and $b$ with weight $\lambda\in(0,1)$, defined by
\[
M_p(a,b;\lambda):=
\begin{cases}
\max\{a,b\}, & \text{if } p=\infty,\\
\bigl((1-\lambda)a^p+\lambda b^p\bigr)^{1/p}, & \text{if } p\in\mathbb{R}\setminus\{0\},\\
a^{1-\lambda}b^\lambda, & \text{if } p=0,\\
\min\{a,b\}, & \text{if } p=-\infty,
\end{cases}
\]
for $a,b>0$, and $M_p(a,b;\lambda):=0$ if $ab=0$. Clearly, if $p=1$, this reduces to the
usual notion of concavity. We often refer to the cases $p=0$ and $p=-\infty$ as
log-concavity and quasi-concavity, respectively.
\end{defi}

Throughout the subsection, we denote
\[
B(x,r):=\{z\in\mathbb{R}^n:|z-x|<r\},\quad x\in\mathbb{R}^n,\ r>0.
\]
For $x=(x_1,\dots,x_n)$ and $y=(y_1,\dots,y_n)\in\mathbb{R}^n$, we write $x\otimes y$ for
the $n\times n$ matrix with entries $(x_i y_j)$, $i,j=1,\dots,n$. For
$m\in\{2,3,\dots\}$, we set
\[
\Lambda_m:=\Bigl\{\lambda=(\lambda_1,\dots,\lambda_m)\in(0,1)^m:
\sum_{i=1}^m\lambda_i=1\Bigr\}.
\]
For $a=(a_1,\dots,a_m)\in(0,\infty)^m$, $\lambda\in\Lambda_m$, and
$p\in[-\infty,\infty]$, we define
\[
{\bf M}_p(a;\lambda):=
\begin{cases}
\bigl[\lambda_1a_1^p+\dots+\lambda_m a_m^p\bigr]^{1/p},
& \text{if } p\neq-\infty,0,\infty,\\
\max\{a_1,\dots,a_m\},
& \text{if } p=\infty,\\
a_1^{\lambda_1}\cdots a_m^{\lambda_m},
& \text{if } p=0,\\
\min\{a_1,\dots,a_m\},
& \text{if } p=-\infty,
\end{cases}
\]
which is the ($\lambda$-weighted ) $p$-mean of $a$; see \cite{IS2014}. Obviously,
${\bf M}_p(a;\lambda)$ generalizes the two-point mean $M_p(a,b;\lambda)$.

In order to study the spatial concavity,  the definition of the $p$-concave $\lambda$-envelope will be introduced, see \cite{INS2020}.

\begin{defi}\label{d2.2}
Let $D$ be a bounded convex domain in $\mathbb{R}^n$, and let $u$ be a nonnegative
function on $[0,\infty)\times\bar D$. Fix $\lambda=(\lambda_1,\dots,\lambda_{n+1})\in\Lambda_{n+1}$ and
$p\in[-\infty,\infty]$. For any $(t,x)\in[0,\infty)\times\bar D$, set
\[
u_{p,\lambda}(t,x)
:=
\sup\Bigl\{
{\bf M}_p\bigl(u(t,y_1),\dots,u(t,y_{n+1});\lambda\bigr):
\{y_i\}_{i=1}^{n+1}\subset\bar D,\
x=\sum_{i=1}^{n+1}\lambda_i y_i
\Bigr\}.
\]
We call $u_{p,\lambda}$ the $p$-concave $\lambda$-envelope of $u$.
\end{defi}
In particular, by taking $y_i\equiv x$ in Definition~\ref{d2.2}, we immediately obtain
\[
u_{p,\lambda}(t,x)\ge u(t,x), \qquad (t,x)\in[0,\infty)\times\bar D.
\]
Consequently, letting $u$ be the solution of \eqref{2.1}, in order to prove that $u_{p,\lambda}$ coincides with $u$ and hence is a viscosity solution of \eqref{2.1},
it suffices to show that
\[
u_{p,\lambda}(t,x)\le u(t,x),\qquad (t,x)\in[0,\infty)\times\bar D.
\]

We now recall the definition of viscosity solutions to \eqref{2.1}.
\begin{defi}\label{d2.3}
(i) An upper semicontinuous function $u$ on $[0,\infty)\times D$ is called a viscosity subsolution of
\eqref{2.1} if, for any $(t_0,x_0)\in[0,\infty)\times D$ and any test function
$\phi\in C^{1,2}([0,\infty)\times D)$ such that
\[
\phi(t_0,x_0)=u(t_0,x_0)
\quad\text{and}\quad
\phi(t,x)\ge u(t,x)\ \text{in a neighborhood of }(t_0,x_0)\ \text{a.s.},
\]
the following inequality holds at $(t_0,x_0)$ almost surely:
\[
\phi(t_0,x_0) \le \phi(0,x_0)+\int_0^{t_0}\Delta\phi(s,x_0)\,ds + \int_0^{t_0}f(s,x_0,\phi(s,x_0))\,ds + \int_0^{t_0}\sigma(s,x_0,\phi(s,x_0))dB(s).
\]

(ii) A lower semicontinuous function $u$ on $[0,\infty)\times D$ is called a viscosity supersolution of
\eqref{2.1} if, for any $(t_0,x_0)\in[0,\infty)\times D$ and any test function
$\phi\in C^{1,2}([0,\infty)\times D)$ such that
\[
\phi(t_0,x_0)=u(t_0,x_0)
\quad\text{and}\quad
\phi(t,x)\le u(t,x)\ \text{in a neighborhood of }(t_0,x_0)\ \text{a.s.},
\]
the following inequality holds at $(t_0,x_0)$ almost surely:
\[
\phi(t_0,x_0) \ge \phi(0,x_0)+\int_0^{t_0}\Delta\phi(s,x_0)\,ds + \int_0^{t_0}f(s,x_0,\phi(s,x_0))\,ds + \int_0^{t_0}\sigma(s,x_0,\phi(s,x_0))dB(s).
\]
\end{defi}

Define
\[
u_p(t,x):=\sup_{\lambda\in\Lambda_{n+1}} u_{p,\lambda}(t,x).
\]
Then, to prove that the strong solution $u$ of \eqref{2.1} is $p$-concave with respect to the spatial variable,
it suffices to show that
\[
u_p(t,x)=u(t,x),\qquad (t,x)\in[0,\infty)\times\bar D.
\]

In this paper, we focus on the log-concavity of solutions to \eqref{2.1}, corresponding to the case $p=0$.
More precisely, our goal is to prove that
\[
u_{0,\lambda}(t,x)=u(t,x),
\qquad (t,x)\in[0,\infty)\times\bar D.
\]

For stochastic partial differential equations of the form \eqref{2.1}, it is well known
that if the reaction term $f$ is monotone and the noise coefficient $\sigma$ is Lipschitz
continuous, then a suitable comparison principle holds.
Inspired by \cite{INS2020}, we obtain the following result.
\begin{theo}\label{t2.3}
Let $D$ be a bounded convex smooth domain in $\mathbb{R}^n$ and $u$ be the strong solution of \eqref{2.1}.
Assume that $f(t,x,u)=f(u)$ satisfies a Lipschitz condition and that
$\sigma(t,x,u)=\sigma u$ with a positive constant $\sigma$.
Suppose that the function $v\mapsto e^{-v}f(e^v)$ is concave.
Then $u_{0,\lambda}$ is a viscosity subsolution of \eqref{2.1}.
Moreover, if the comparison principle holds for \eqref{2.1}, then
\[
u(t,x)=u_{0,\lambda}(t,x),
\qquad (t,x)\in[0,T]\times\bar D,
\]
for any $T>0$, and $u$ is log-concave with respect to the spatial variable.
\end{theo}

\begin{proof} It follows from the classical theory \cite{Cb2007,Da1992} that under the assumptions of Theorem \ref{t2.3}, the \eqref{2.1} admits a unique strong solution.
Fix $(t_*,x_*)\in(0,T]\times D$ and $\lambda=(\lambda_1,\dots,\lambda_{n+1})\in\Lambda_{n+1}$.
By \cite[Theorem~4.4]{KY2020}, we have
\[
u(t_*,x_*)>0
\quad\text{almost surely.}
\]
The set
\[
\Bigl\{\{y_i\}_{i=1}^{n+1}\subset\bar D:\;
x_*=\sum_{i=1}^{n+1}\lambda_i y_i\Bigr\}
\]
is closed and bounded in $\bar D$, and the mapping
\[
(y_1,\dots,y_{n+1})
\longmapsto
{\bf M}_0\bigl(u(t_*,y_1),\dots,u(t_*,y_{n+1});\lambda\bigr)
\]
is continuous.
Using the boundary condition $u(t,\cdot)|_{\partial D}=0$, we can choose
$\{x_i\}_{i=1}^{n+1}\subset D$ such that
\[
x_*=\sum_{i=1}^{n+1}\lambda_i x_i
\]
and
\[
u_{0,\lambda}(t_*,x_*)
=
{\bf M}_0\bigl(u(t_*,x_1),\dots,u(t_*,x_{n+1});\lambda\bigr)
=
\prod_{i=1}^{n+1}u(t_*,x_i)^{\lambda_i}.
\]
By the Lagrange multiplier method, we obtain
\[
\frac{\nabla u(t_*,x_1)}{u(t_*,x_1)}
=
\cdots
=
\frac{\nabla u(t_*,x_{n+1})}{u(t_*,x_{n+1})}
=:\Theta.
\]

Set
\[
u_*:=u_{0,\lambda}(t_*,x_*),\qquad
u_i:=u(t_*,x_i),\qquad
y_i(x):=x_i+(x-x_*).
\]
Then $x=\sum_{i=1}^{n+1}\lambda_i y_i(x)$.
Define
\[
\phi(t,x)
:=
{\bf M}_0\bigl(u(t,y_1(x)),\dots,u(t,y_{n+1}(x));\lambda\bigr)
=
\prod_{i=1}^{n+1}u(t,y_i(x))^{\lambda_i}.
\]
The function $\phi$ is $C^2$ with respect to the spatial variable almost surely and
satisfies $\phi(t_*,x_*)=u_*$.
Moreover, by the definition of $u_{0,\lambda}$,
\[
u_{0,\lambda}(t,x)\ge \phi(t,x)
\]
in a neighborhood of $(t_*,x_*)$ almost surely.

Since $f$ satisfies the Lipschitz condition, the classical comparison principle
can be applied to \eqref{2.1}, see \cite{LD2015}.
It therefore suffices to show that
\begin{equation}\label{2.6}
\phi(t_*,x_*) \le \phi(0,x_*)+\int_0^{t_*}\Delta\phi(s,x_*)\,ds + \int_0^{t_*}f(\phi(s,x_*))\,ds + \int_0^{t_*}\sigma(s,x_*,\phi(s,x_*))dB(s).
\end{equation}

A direct computation yields
\bess
d\phi(t_*,x_*)
&=&
\phi(t_*,x_*)\sum_{i=1}^{n+1}\lambda_i
\frac{du(t_*,x_i)}{u_i},\\[0.3em]
\nabla\phi(t_*,x_*)
&=&
\phi(t_*,x_*)\sum_{i=1}^{n+1}\lambda_i
\frac{\nabla u(t_*,x_i)}{u_i}
=\Theta,\\[0.3em]
\Delta\phi(t_*,x_*)
&=&
\phi(t_*,x_*)\sum_{i=1}^{n+1}\lambda_i
\frac{\Delta u(t_*,x_i)}{u_i},
\eess
where we used the identity
\[
\sum_{i=1}^{n+1}\lambda_i(\lambda_i-1)
+
\sum_{\substack{1\le i,j\le n+1\\ i\neq j}}
\lambda_i\lambda_j
=0.
\]
Consequently,
\begin{equation}\label{2.7}
\begin{aligned}
&\phi(t_*,x_*)-\phi(0,x_*)-\int_0^{t_*}\Delta\phi(s,x_*)\,ds - \int_0^{t_*}f(\phi(s,x_*))\,ds - \int_0^{t_*}\sigma(s,x_*,\phi(s,x_*))dB(s)\\
&\qquad=\int_0^{t_*}
\Biggl[
\phi(s,x_*)\sum_{i=1}^{n+1}\lambda_i
\frac{f(u_i)}{u_i}
-
f(\phi(s,x_*))
\Biggr]ds .
\end{aligned}
\end{equation}
Let $v_i=\log u_i$.
By the assumed concavity of $v\mapsto e^{-v}f(e^v)$, we have
\bess
\sum_{i=1}^{n+1}\lambda_i\frac{f(u_i)}{u_i}
=
\sum_{i=1}^{n+1}\lambda_i e^{-v_i}f(e^{v_i})
\le
e^{-\sum_{i=1}^{n+1}\lambda_i v_i}
f\!\left(e^{\sum_{i=1}^{n+1}\lambda_i v_i}\right)
=
\phi(s,x_*)^{-1}f(\phi(s,x_*)).
\eess
Substituting this inequality into \eqref{2.7} yields
\[
\phi(t_*,x_*)-\phi(0,x_*)-\int_0^{t_*}\Delta\phi(s,x_*)\,ds - \int_0^{t_*}f(\phi(s,x_*))\,ds - \int_0^{t_*}\sigma(s,x_*,\phi(s,x_*))dB(s)
\le 0
\quad\text{a.s.}
\]
This proves \eqref{2.6} and completes the proof.
\end{proof}

\begin{remark}\label{r2.1}

(1)
The assumption $\sigma(u)=\sigma u$ avoids difficulties caused by the lack of order
preservation for stochastic integrals. More precisely, even if random variables $X>Y$ almost surely, one cannot in general conclude that
\[
XdB(t)\ge YdB(t)\quad\text{a.s.}
\]
In this case $\sigma(u)=\sigma u$, the noise term disappears after a logarithmic transformation, allowing us to deduce almost sure spatial log-concavity.

(2)
By \cite{KY2020}, solutions of \eqref{2.1} remain positive almost surely.
Therefore, the heat energy
\[
H(t)=\int_D \mathbb{E}[u(t,x)]dx
\]
is well defined and is log-concave with respect to time.

(3)
An example of a function $f$ satisfying the structural condition in
Theorem~\ref{t2.3} is
\[
f(u)=-\sqrt{1+u^2},\qquad u>0.
\]

(4) If $\sigma=0$, Theorem~\ref{t2.3} also applies to the deterministic equation \eqref{2.1},
yielding preservation of spatial log-concavity in the deterministic setting.
\end{remark}
\section{Nonlinear Multiplicative Noise}\setcounter{equation}{0}

In this section, we study the impact of noise of geometry property of parabolic equation on bounded domain. For simplicity, we only consider scalar case. Consider
\bes
\left\{
\begin{array}{llll}
d u(t,x)
=
\bigl[\Delta u(t,x)+f(u(t,x),\nabla u(t,x))\bigr]dt
+\sigma(u(t,x),\nabla u(t,x))dB(t),
& x\in D,\ t>0,\\
u(t,x)=0,
& x\in\partial D,\ t>0,\\
u(0,x)=u_0(x),
& x\in D,
\end{array}
\right.
\label{3.1}
\ees
where $B(t)$ is one-dimensional Brownian motion similar to that in (\ref{2.1}). It follows from section 2 that it is impossible to consider the geometry property of equation \eqref{3.1} almost surely. Now, we are ready to consider the moment of the solution of \eqref{3.1}.

By using It\^{o} formula, we have, for a function $\psi$,
  \bess
d\psi(u)&=&\psi'(u)du+\frac{1}{2}\psi''(u)du\cdot du\\
&=&\left[\psi'(u)\Delta u+\psi'(u)f(u,\nabla u)+\frac{1}{2}\psi''(u)\sigma^2(u)\right]dt+\psi'(u)\sigma(u,\nabla u)dB(t)\\
&=&\left[\Delta \psi(u)-\psi''(u)|\nabla u|^2+\psi'(u)f(u,\nabla u)+\frac{1}{2}\psi''(u)\sigma^2(u,\nabla u)\right]dt+\psi'(u)\sigma(u,\nabla u)dB(t).
  \eess
Obviously, if, for any $(t,x)\in(0,\infty)\times D$, it holds that
  \bess
-\psi''(u)|\nabla u|^2+\psi'(u)f(u,\nabla u)+\frac{1}{2}\psi''(u)\sigma^2(u,\nabla u)=0,
  \eess
then letting $v=\psi(u)$ and denote $u=\psi^{-1}(v)$ (we assume the inverse function exists), we get
\bes
\left\{
\begin{array}{llll}
d v(t,x)
=\Delta v(t,x) dt
+\psi'(\psi^{-1}(v(t,x)))\sigma(\psi^{-1}(v(t,x)),\nabla\psi^{-1}(v(t,x)))dB(t),
\\
v(t,x)=\psi(0),
\\
v(0,x)=\psi(u_0(x)).
\end{array}
\right.
\label{3.2}
\ees
Moreover, if we take $w=\mathbb{E}[v]-\psi(0)$, then $w$ satisfies
\bes
\left\{
\begin{array}{llll}
\partial_t w(t,x)
=\Delta w(t,x),
& x\in D,\ t>0,\\
w(t,x)=0,
& x\in\partial D,\ t>0,\\
w(0,x)=\psi(u_0(x))-\psi(0),
& x\in D.
\end{array}
\right.
\label{3.3}
\ees
Ishige et al. \cite{IST2024} considered the F-concavity of Dirichlet flow \eqref{3.3}.  In order to write it clearly, following \cite{IST2024,I2024}, we introduce the notion of $F$-concavity.
\begin{defi}\label{d3.1} Let $a\in(0,\infty]$, and set $I=[0,a)$ and int$I=(0,a)$. A function $F:I\to[-\infty,\infty)$ is said admissible on $I$ if $F\in C(\rm{int}\ I)$, $F$ is strictly increasing on $I$, and $F(0)=-\infty$. Let $\mathcal{A}(I)$ be the set of $I$-valued functions in $\mathbb{R}^n$. Let $F$ be admissible on $I$. For any $f\in\mathcal{A}(I)$, we say that $f$ is $F$-concave in $D$ if
  \bess
F(f(\lambda x+(1-\lambda)y))\geq \lambda F(f(x))+(1-\lambda)F(f(y))
  \eess
for $x,y\in D$ and $\lambda\in(0,1)$. We denote by $\mathcal{C}[F]$ the set of $F$-concave functions in $\mathbb{R}^n$.
\end{defi}
Let 
  \bess
h(z):=(4\pi)^{-\frac{1}{2}}\int_0^\infty e^{-\frac{|z-w|^2}{4}}dw\ \ \ {\rm for}\ z\in\mathbb{R}.
 \eess
Then the function $h$ is smooth in $\mathbb{R}$, $\lim_{z\to-\infty}h(z)=0$, $\lim_{z\to\infty}h(z)=1$, and $h'>0$ in $\mathbb{R}$. Thus the inverse function of $h$ exists, and denote by $H$. Define an admissible function $H_a$ on $[0,a)$ by 
  \bess
H_a(r):=\left\{\begin{array}{llll}
H(\frac{r}{a})\ \ \ & {\rm for}\ r\in(0,a)\ {\rm if}\ a>0,\\
\log r\  \ & {\rm for}\ r\in(0,a)\ {\rm if}\ a=\infty,\\
-\infty\  \ & {\rm for}\ r=0\ {\rm and}\ a\in(0,\infty].
  \end{array}\right.\eess
"The F-concavity is preserved by Dirichlet heat flow (DHF)" means that the solution of \eqref{3.3} is  F-concavity if the initial data is F-concavity. As for \eqref{3.3}, we have the following result. 
  \begin{theo}\label{t3.1} Let $I=[0,a)$ with $a\in(0,\infty]$ and $D$ be a convex domain in $\mathbb{R}^n$ with $n\geq1$. Then the following properties hold:

(i) $H_a$-concavity is the strongest $F$-concavity preserved by DHF in $\mathcal{A}_D$;

(ii) Log-concavity is the weakest $F$-concavity preserved by DHF in $\mathcal{A}_D$ for $n\geq2$; when $n=1$, the same is true under $C^2$-regularity assumption on $F$;

(iii) When $n\geq2$, if starting with an initial datum which is not log-concave, even quasi-concavity can be immediately broken by the corresponding solution of problem \eqref{3.3}, so losing every possible reminiscence of concavity.
  \end{theo}

The proof of Theorem \ref{t3.1} is the exact same as that of Theorem 1.5 of \cite{IST2024}. 

Now, we want to show the impact of noise on F-concavity. It follows from (ii) of of Theorem \ref{t3.1} that if $\psi(u_0(x))-\psi(0)$ is log-concavity, then solution $w$ is also log-concavity. Consider 
\bes
\left\{
\begin{array}{llll}
d u(t,x)
=
\bigl[\Delta u(t,x)-u(t,x)-2\sqrt{2}|\nabla u(t,x)|)\bigr]dt\\
\qquad\quad+[u(t,x)+\sqrt{2}|\nabla u(t,x)|]dB(t),
\ \  x\in D,\ t>0,\\
u(t,x)=0,
 \ \  x\in\partial D,\ t>0,\\
u(0,x)=u_0(x),
\ \  x\in D,
\end{array}
\right.
\label{3.4}
\ees
Let $\psi(x)=x^2$. Direct calculations show that $w=\mathbb{E}[\psi(u)]-\psi(0)$ satisfies
\bess
\left\{
\begin{array}{llll}
\partial_t w(t,x)
=\Delta w(t,x),
& x\in D,\ t>0,\\
w(t,x)=0,
& x\in\partial D,\ t>0,\\
w(0,x)=u_0(x)^2,
& x\in D.
\end{array}
\right.
\eess
Obviously, it follows from Theorem \ref{t3.1} that if $u_0^2$ satisfies Log-concavity, then the solution $\mathbb{E}[u^2(t,x)]$ is also Log-concavity in spatial variable. But if $\sigma=0$, then $w(t,x)=u^2(t,x)$ satisfies 
\bes
\left\{
\begin{array}{llll}
\partial_t w(t,x)
=\Delta w(t,x)-\frac{|\nabla w(t,x)|^2}{2w(t,x)}-2w(t,x)-\sqrt{2}|\nabla w(t,x)|,
& x\in D,\ t>0,\\
w(t,x)=0,
& x\in\partial D,\ t>0,\\
w(0,x)=u_0(x)^2,
& x\in D.
\end{array}
\right.
\label{3.5}
\ees
Set
  \bess
F(w,\nabla w,\Delta w)=\Delta w--\frac{|\nabla w|^2}{2w}-2w-\sqrt{2}|\nabla w|.
 \eess
Then it follows from \cite[Theorem 4.3]{INS2020} that the solution of \eqref{3.5} is Log-concavity under the condition that, for any fixed $\theta\in\mathbb{R}^n$,
  \bess
e^{-w}F(e^w,e^w\theta,e^wA)
  \eess
is convex with respect to $(w,A)\in \mathbb{R}\times S^n$, where $S^n$ denotes the space
of $n\times n$ real symmetric matrices. Unfortunately, 
  \bess
e^{-w}F(e^w,e^w\theta,e^wA)=A-e^{-w}|\theta|-2-\sqrt{2}|\theta|,
  \eess
which is not convex because the term $-e^{-w}$ is concave. Therefore, the noise can help the solution keep the same property as initial data.

\section{Generalized Pr\'{e}kopa-Leindler inequality}\label{Sec4}
\setcounter{equation}{0}
In Sections~\ref{Sec2}, we established preservation results for spatial concavity and log-concavity of solutions to stochastic parabolic equations on bounded convex domains.
It is well known that geometric properties of parabolic flows are closely connected with functional inequalities of Brunn--Minkowski type.
In the whole space $\mathbb{R}^n$, the Pr\'ekopa-Leindler inequality can be interpreted as an integral manifestation of the log-concavity preserving property of the heat semigroup.

Motivated by this perspective, we show in this section that the log-concavity preserving structure of the Dirichlet heat flow on bounded convex domains yields a generalized Pr\'ekopa-Leindler inequality involving the first Dirichlet eigenfunction.
This provides a functional counterpart to the pointwise geometric preservation results obtained earlier.

We first recall the Borell-Brascamp-Lieb inequality. Let $n\in\mathbb{N}$, $\lambda\in(0,1)$ and $\alpha\in[-1/n,\infty]$.
Assume that nonnegative integrable functions $f,g$ and $h$ in $\mathbb{R}^n$ satisfy
\[
h((1-\lambda)y+\lambda z)\geq M_\alpha(f(y),g(z);\lambda)
\]
for almost all $y,z\in\mathbb{R}^n$. Then
\begin{equation}
\int_{\mathbb{R}^n}h(x)dx
\geq
M_{\frac{\alpha}{1+n\alpha}}
\left(
\int_{\mathbb{R}^n}f(x)dx,
\int_{\mathbb{R}^n}g(x)dx;\lambda
\right).
\label{4.1}
\end{equation}

When $\alpha=0$, this reduces to the classical Pr\'ekopa-Leindler inequality:
\begin{equation}
\int_{\mathbb{R}^n}h(x)dx
\geq
\left(\int_{\mathbb{R}^n}f(x)dx\right)^{1-\lambda}
\left(\int_{\mathbb{R}^n}g(x)dx\right)^{\lambda}.
\label{4.2}
\end{equation}

\medskip
Consider the Dirichlet heat equation
\begin{equation}
\left\{
\begin{aligned}
&\partial_t u(t,x)= \Delta u(t,x), && x\in D,\ t>0,\\
&u(t,x)=0, && x\in\partial D,\ t>0,\\
&u(0,x)=u_0(x), && x\in D,
\end{aligned}
\right.
\label{4.3}
\end{equation}
where $D\subset\mathbb{R}^n$ is bounded convex smooth. The solution admits the spectral representation
\begin{equation}
u(t,x)=\sum_{n=1}^\infty c_n e^{-\lambda_n t}\phi_n(x),
\label{4.4}
\end{equation}
where $c_n=\int_D u_0(x)\phi_n(x)dx$ and $(\lambda_n,\phi_n)$ are the Dirichlet eigenpairs of $-\Delta$:
\[
-\Delta\phi_n=\lambda_n\phi_n,\quad \phi_n|_{\partial D}=0,
\]
with
\[
0<\lambda_1<\lambda_2\leq\lambda_3\leq\cdots.
\]

\begin{theo}\label{t4.1}
Let $D$ be bounded convex smooth.
Assume $\psi_0,\psi_1,\psi_\lambda\in C(D)$ are nonnegative and satisfy
\begin{equation}
\psi_\lambda(\lambda x+(1-\lambda)y)
\geq
\psi_0(x)^\lambda \psi_1(y)^{1-\lambda}.
\label{4.5}
\end{equation}
Let $u_i$ solve \eqref{4.3} with initial data $\psi_i$. Then
\begin{equation}
u_\lambda(t,\lambda x+(1-\lambda)y)
\geq
u_0(t,x)^\lambda u_1(t,y)^{1-\lambda}
\label{4.6}
\end{equation}
for all $x,y\in D$ and $t>0$.
\end{theo}

\begin{proof}
Define
\[
u_{0,\lambda}(t,x)
=
\sup\left\{
M_0(u_0(t,y),u_1(t,z);\lambda):
y,z\in\bar D,\
x=\lambda y+(1-\lambda)z
\right\}.
\]
By the same viscosity-envelope argument used in the proof of Theorem~\ref{t2.3}, replacing $u(t,x_i)$ by $u_i(t,x_i)$, one verifies that $u_{0,\lambda}$ is a viscosity solution of \eqref{4.3}.

By uniqueness of solutions to the Dirichlet heat equation, $u_{0,\lambda}$ coincides with the solution of \eqref{4.3} with initial data
\[
\psi_{0,\lambda}
=
\sup\left\{
M_0(\psi_0(y),\psi_1(z);\lambda):
y,z\in\bar D,\
x=\lambda y+(1-\lambda)z
\right\}.
\]

Condition \eqref{4.5} implies $\psi_\lambda\geq \psi_{0,\lambda}$.
By the comparison principle,
\[
u_\lambda(t,\lambda x+(1-\lambda)y)
\geq
u_{0,\lambda}(t,x).
\]
By definition of $u_{0,\lambda}$,
\[
u_{0,\lambda}(t,x)
\geq
u_0(t,x)^\lambda u_1(t,y)^{1-\lambda}.
\]
The proof is complete.
\end{proof}

\medskip
Using the asymptotic dominance of the first eigenfunction in \eqref{4.4}, we now derive a weighted Pr\'ekopa-Leindler inequality on $D$.

\begin{theo}\label{t4.2}
Assume nonnegative integrable functions $f,g,h$ on convex domain $D$ satisfy
\[
h((1-\lambda)y+\lambda z)
\geq
f(y)^\lambda g(z)^{1-\lambda}
\]
for almost all $y,z\in D$. Then
\begin{equation}
\int_D h(x)\phi_1(x)dx
\geq
\left(
\int_D f(x)\phi_1(x)dx
\right)^\lambda
\left(
\int_D g(x)\phi_1(x)dx
\right)^{1-\lambda}.
\label{4.7}
\end{equation}
\end{theo}

\begin{proof}
First assume $f,g,h\in C(D)$.
Let $u_\lambda,u_0,u_1$ solve \eqref{4.3} with initial data $h,f,g$. By Theorem~\ref{t4.1},
\[
\sum_{n=1}^\infty c_{n,\lambda}e^{-\lambda_n t}\phi_n(x)
\geq
\left(\sum_{n=1}^\infty c_{n,0}e^{-\lambda_n t}\phi_n(x)\right)^\lambda
\left(\sum_{n=1}^\infty c_{n,1}e^{-\lambda_n t}\phi_n(x)\right)^{1-\lambda}.
\]

Write
\[
u_i(t,x)
=
e^{-\lambda_1 t}
\left(
c_{1,i}\phi_1(x)
+
\sum_{n=2}^\infty
c_{n,i}e^{-(\lambda_n-\lambda_1)t}\phi_n(x)
\right).
\]

Since $\lambda_n>\lambda_1$ for $n\ge2$, the higher modes decay exponentially.
Dividing by $e^{-\lambda_1 t}$ and letting $t\to\infty$, we obtain
\[
c_{1,\lambda}\phi_1(x)
\geq
(c_{1,0}\phi_1(x))^\lambda
(c_{1,1}\phi_1(x))^{1-\lambda}.
\]

Recalling
\[
c_{1,i}
=
\int_D f_i(x)\phi_1(x)dx,
\]
we obtain \eqref{4.7}.
The extension to $L^1$ data follows by approximation.
\end{proof}

We now return to the stochastic setting considered in Section~\ref{Sec2}.
Consider
\begin{equation}
\left\{
\begin{aligned}
&du(t,x)= \Delta u(t,x)dt + \sigma u(t,x)dB(t),\\
&u(t,x)|_{\partial D}=0.
\end{aligned}
\right.
\label{4.8}
\end{equation}

Define
\[
v(t,x)=e^{-\sigma B(t)}u(t,x).
\]
Then $v$ satisfies
\[
\partial_t v(t,x)
=
\Delta v(t,x) -\frac{\sigma^2}{2}v(t,x).
\]

Let $w(t,x)=e^{\frac{\sigma^2}{2}t}v(t,x)$.
Then $w$ solves the heat equation \eqref{4.3} almost surely. Applying Theorem~\ref{t4.2} pathwise yields:
\begin{theo}\label{t4.3}
Let $f,g,h$ be nonnegative integrable random functions satisfying
\[
h((1-\lambda)y+\lambda z)
\geq
f(y)^\lambda g(z)^{1-\lambda}
\]
almost surely.
Then
\[
\int_D h(x)\phi_1(x)dx
\geq
\left(
\int_D f(x)\phi_1(x)dx
\right)^\lambda
\left(
\int_D g(x)\phi_1(x)dx
\right)^{1-\lambda}
\]
almost surely.
\end{theo}
In the end of this section, we give another method to prove the inequality \eqref{4.7}. First, we use the zero extension method to obtain the Pr\'ekopa-Leindler inequality on bounded domain.  Assume the assumptions of Theorem \ref{t4.2} hold, then define 
  \bess
w_s(x)=\left\{\begin{array}{lll}
w(x), \ \ \ & if\ x\in D,\\
0,\ \ \ \ & if x\in\mathbb{R}^n\setminus D,
\end{array}\right.\ \ w=f\ {\rm or}\ g\ {\rm or}\ h.
  \eess
Then $f_s,g_s,h_s$ satisfy the assumption of the Pr\'ekopa-Leindler inequality in whole space. That is to say, $f_s,g_s,h_s$ are nonnegative integrable functions satisfying
\[
h((1-\lambda)y+\lambda z)\geq f(y)^\lambda g(z)^{1-\lambda}
\]
for almost all $y,z\in\mathbb{R}^n$.
Conseqently, we have 
\[
\int_{\mathbb{R}^n} h_s(x)dx
\geq
\left(
\int_{\mathbb{R}^n} f_s(x)dx
\right)^\lambda
\left(
\int_{\mathbb{R}^n} g_x(x)dx
\right)^{1-\lambda}.
\]
which implies 
\[
\int_D h(x)dx
\geq
\left(
\int_D f(x)dx
\right)^\lambda
\left(
\int_D g(xdx
\right)^{1-\lambda}.
\]
Note that the first eigenvalue function $\phi_1$ is Log-concavity on convex domain \cite{K198301,K198302}, i.e.,
   \bess
\phi_1((1-\lambda)y+\lambda z)\geq \phi_1(y)^\lambda \phi_1(z)^{1-\lambda}
  \eess
for almost all $y,z\in D$. Combining the assumptions on $f,g,h$, we get 
 \bess
(h\phi_1)((1-\lambda)y+\lambda z)\geq (f\phi_1)(y)^\lambda (g\phi_1)(z)^{1-\lambda}
  \eess
for almost all $y,z\in D$. Applying the Pr\'ekopa-Leindler inequality on bounded domain, we get \eqref{4.7}. However, in this paper, we give a new methd to get the inequality \eqref{4.7}.

\section*{Use of AI tools declaration}
We did not use any AI tools to assist in the derivation and writing of this paper.

\section*{Acknowledgements} This work was supported by the NSF of China grants 12171247.

\section*{Data Availability Statements} Data sharing is not applicable to this article as no new data were created or analyzed in this study.

\section*{Conflicts of Interest} The authors declare that they have no competing interests.

 \end{document}